\documentclass{amsart}
\usepackage{mathrsfs}
\usepackage[all]{xy}
\usepackage{amssymb}

\newcommand{\C}{\mathbb{C}}
\newcommand{\D}{\mathbb{D}}
\newcommand{\F}{\mathbb{F}}
\newcommand{\Gm}{\mathbb{G}_m}
\newcommand{\Q}{\mathbb{Q}}
\newcommand{\Qlb}{\bar{\mathbb{Q}}_\ell}
\newcommand{\Z}{\mathbb{Z}}

\newcommand{\cA}{\mathscr{A}}
\newcommand{\cB}{\mathscr{B}}

\newcommand{\cF}{\mathcal{F}}

\newcommand{\cN}{\mathcal{N}}

\newcommand{\tcN}{{\tilde{\mathcal{N}}}}
\newcommand{\cO}{\mathcal{O}}

\newcommand{\fg}{\mathfrak{g}}

\newcommand{\ft}{\mathfrak{t}}
\newcommand{\fu}{\mathfrak{u}}
\newcommand{\fub}{\bar{\mathfrak{u}}}

\newcommand{\bb}{\mathbf{1}}
\newcommand{\bA}{\mathbf{A}}

\newcommand{\bH}{\mathrm{H}}
\newcommand{\bmH}{\breve{\mathrm{H}}}

\newcommand{\IC}{\mathrm{IC}}

\newcommand{\Dbc}{\mathrm{D}^{\mathrm{b}}_{\mathrm{m}}}
\newcommand{\DbcGm}{\mathrm{D}^{\mathrm{b}}_{\mathrm{m},\Gm}}
\newcommand{\Spr}{\mathbf{Spr}}
\newcommand{\Dspr}{\mathrm{D}_{\Spr}}
\newcommand{\Dsprp}{\mathrm{D}_{\Spr}^\perp}

\newcommand{\pt}{\mathrm{pt}}

\newcommand{\Forg}{\mathbb{U}}

\newcommand{\PN}{\mathrm{P}_G(\cN)}
\newcommand{\Sh}{\mathrm{Sh}_G(\cB)}
\newcommand{\Shpt}{\mathrm{Sh}(\pt)}

\newcommand{\Rep}{\mathrm{Rep}}
\newcommand{\Lie}{\mathrm{Lie}}

\newtheorem{thm*}{Theorem}
\numberwithin{equation}{section}
\newtheorem{thm}{Theorem}[section]
\newtheorem{lem}[thm]{Lemma}
\newtheorem{prop}[thm]{Proposition}
\newtheorem{cor}[thm]{Corollary}
\theoremstyle{definition}

\theoremstyle{remark}
\newtheorem{rmk}[thm]{Remark}

\DeclareMathOperator{\Hom}{Hom}
\DeclareMathOperator{\End}{End}
\DeclareMathOperator{\Aut}{Aut}

\DeclareMathOperator{\Irr}{Irr}

\DeclareMathOperator{\chr}{char}
\DeclareMathOperator{\uHom}{\underline{Hom}}
\DeclareMathOperator{\cRHom}{\mathit{R}\mathcal{H}\mathit{om}}
\DeclareMathOperator{\RHom}{\mathit{R}\mathrm{Hom}}
\DeclareMathOperator{\uEnd}{\underline{End}}

\DeclareMathOperator{\Spec}{Spec}

\newcommand{\simto}{\mathrel{\overset{\sim}{\to}}}

\title{Green functions via hyperbolic localization}
\author{Pramod N. Achar}
\address{Department of Mathematics, Louisiana State University, Baton Rouge, LA 70803, USA}
\email{pramod@math.lsu.edu}
\thanks{The author was partly supported by National Security
Agency grant H98230-09-1-0024.}

\begin{document}

\begin{abstract}
Let $G$ be a reductive algebraic group, with nilpotent cone $\cN$ and flag variety $\cB$.  We construct an exact functor from perverse sheaves on $\cN$ to locally constant sheaves on $\cB$, and we use it to study Ext-groups and stalks of simple perverse sheaves on $\cN$ in terms of the cohomology of $\cB$.
\end{abstract}

\maketitle

\section{Introduction}
\label{sect:intro}

Let $G$ be a connected reductive algebraic group over an algebraically closed field $\Bbbk$ of good characteristic.  Let $\cN$ denote the nilpotent cone in its Lie algebra $\fg$, and let $W$ denote its Weyl group.  An explicit description of the stalks of simple perverse sheaves on $\cN$ has been given by Lusztig~\cite{lus:cs5}, building on earlier ideas of Shoji~\cite{sho:gpF4,sho:gpcg}.  For most such perverse sheaves (those appearing in the Springer correspondence), this description involves the representation theory of $W$, and specifically its \emph{coinvariant algebra}.  The coinvariant algebra of $W$ is also isomorphic to the cohomology ring $\bH^\bullet(\cB)$ of the flag variety $\cB$.  However, that cohomology ring does not appear in the proofs in~\cite{lus:cs5}, which rely instead on orthogonality properties of character sheaves coming from the geometry of semisimple classes.

The present paper is an attempt to understand Lusztig's results directly in terms of the geometry of $\cB$.  Consider the diagram
\begin{equation}\label{eqn:maindiag}
\xymatrix{
\cN & \tcN \ar[l]_{\mu} \ar[r]^{\pi} & \cB,}
\end{equation}
where $\tcN$ denotes the cotangent bundle of $\cB$, $\pi$ is the natural projection map, and $\mu$ is the Springer resolution.  We study the functor
\begin{equation}\label{eqn:phi-defn}
\Phi = \pi_!\mu^*: \Dbc(\cN) \to \Dbc(\cB),
\end{equation}
where $\Dbc(X)$ denotes the category of mixed complexes of $\Qlb$-sheaves (for some $\ell \ne \chr \Bbbk$) that are constructible with respect to the $G$-orbits on $X$. 

The main results, proved in Sections~\ref{sect:exact}--\ref{sect:weyl}, are summarized in the statement below.  This statement involves the following categories:  $\PN \subset \Dbc(\cN)$ is the abelian category of perverse sheaves; $\Spr \subset \PN$ is the Serre subcategory containing the simple perverse sheaves appearing in the Springer correspondence; and $\Sh \subset \Dbc(\cB)$ is the abelian category of locally constant sheaves.

\begin{thm}\label{thm:main}
The functor $\Phi$ restricts to give an exact functor of abelian categories $\Phi|_{\PN}: \PN \to \Sh$.  Moreover, for $F, F' \in \Spr$, the objects $\Phi(F)$ and $\Phi(F')$ and the vector space $\uHom^i(\Phi(F),\Phi(F'))$ each carry a natural action of the Weyl group $W$, and $\Phi$ induces an isomorphism
\begin{equation}\label{eqn:ext-winv}
\uHom_{\Dbc(\cN)}^i(F,F') \simto \uHom_{\Dbc(\cB)}^i(\Phi(F),\Phi(F'))^W.
\end{equation}
\end{thm}

Here, the notation ``$\uHom$'' denotes a Hom-group equipped with an action of Frobenius; along the way to the theorem above, we show that $\uHom^i_{\Dbc(\cN)}(F,F')$ is pure.  However, weights and purity are not used in any essential way; the main results are also valid in the unmixed setting.

The $W$-action on $\Phi(F)$ induces one on the space of global sections $\Gamma(\Phi(F))$, and the composition $\Gamma \circ \Phi: \PN \to \Rep(W)$ turns out to be an equivalence of categories that may be regarded as a categorical version of the Springer correspondence.
 On the other hand, the $W$-action on $\uHom^i(\Phi(F),\Phi(F'))$ is a generalization of the usual action of $W$ on $\bH^\bullet(\cB)$.  Indeed,~\eqref{eqn:ext-winv} can be used together with known formulas for the \emph{fake degrees} of $W$ to carry out explicit calculations of Ext-groups.

As an application, in Sections~\ref{sect:orth}--\ref{sect:green}, we use Theorem~\ref{thm:main} to give new proofs of two results from~\cite{lus:cs5}: a decomposition of $\Dbc(\cN)$ into orthogonal subcategories, and the algorithmic description of stalks of perverse sheaves on $\cN$ mentioned above.

\subsection*{Acknowledgments}
The author is grateful to J.M.~Douglass, A.~Henderson, and D.~Treumann for helpful conversations.

\section{Preliminaries}
\label{sect:prelim}

\subsection{General conventions}
\label{subsect:var}

Throughout the paper, $G$ and all related varieties will be assumed to be defined over the algebraic closure $\Bbbk$ of a finite field $\F_q$ and equipped with an $\F_q$-rational structure.  For any $G$-variety $X$, $\Dbc(X)$ will denote the category of mixed \'etale $\Qlb$-complexes that are constructible with respect to some $G$-stable stratification.  Let $\bb_X$ denote the constant sheaf with value $\Qlb$, and let $\pt = \Spec \Bbbk$.  Let $\omega_X = a^!\bb_\pt$ (where $a: X \to \pt$ is the constant map) denote the dualizing complex, and let $\D = \cRHom(\cdot, \omega_X)$ denote the Verdier duality functor.

For $F, F'\in  \Dbc(X)$, we let $\uHom^i(F,F') = H^i(a_*\cRHom(F,F'))$.  This is a mixed $\Qlb$-sheaf on a point, i.e., a $\Qlb$-vector space equipped with an action of Frobenius.  Forgetting that action yields the $\Qlb$-vector space of morphisms $F \to F'[i]$ over $\Bbbk$.


Nearly all results (those not explicitly involving purity) are also valid in the setting of unmixed sheaves over an arbitrary algebraically closed field $\Bbbk$ of good characteristic, or for $\C$-sheaves in the classical topology when $\Bbbk = \C$.

\subsection{Further notation}
\label{subsect:fnot}

Let $d = \dim \cB$.  We will frequently encounter shifts and Tate twists related to $\dim \cB$, so we adopt the following shorthand notation: if $F$ is any sheaf, morphism, or functor, we put
\[
F^\flat = F[2d](d)
\qquad\text{and}\qquad
F^\sharp = F[-2d](-d).
\]
Because $\cB$ is a smooth variety of dimension $d$, we have $\omega_\cB \simeq \bb_\cB^\flat$, and because $\pi$ is a smooth map of relative dimension $d$, we also have $\pi^! \simeq (\pi^*)^\flat$ and $\omega_\tcN^\sharp \simeq \bb_\tcN^\flat$.

Throught the paper, $W$ will denote the \emph{universal Weyl group} of $G$, cf.~\cite[\S 3.1]{cg:rtcg}.  This group does not depend on the choice of a maximal torus or a Borel subgroup, and we do not fix any such choice in this paper.  Let $\Irr(W)$ denote the set of isomorphism classes of irreducible representations of $W$ on $\Qlb$-vector spaces.  For each $\chi \in \Irr(W)$, choose a representative $V_\chi$.  We will sometimes regard $V_\chi$ as a pure object of weight $0$ in $\Dbc(\pt)$, by letting the Frobenius act on it as the identity.  It is known that
\begin{equation}\label{eqn:contragr}
V_\chi \simeq V_\chi^*
\end{equation}
for all $\chi \in \Irr(W)$.  However, this isomorphism is not canonical.

It is well known that $W$ acts naturally on the cohomology ring $\bH^\bullet(\cB)$, and that under this action, $\bH^\bullet(\cB)$ can be identified with the \emph{coinvariant algebra} of $W$.  Let
\[
\kappa: W \to \Aut(\bH^\bullet(\cB))
\]
denote this action.

Let $Z$ denote the Steinberg variety $Z = \tcN \times_\cN \tcN$.  Finally, let $\iota: \cB \to \tcN$ be the inclusion of the zero section, and let $i_0: \pt \to \cN$ denote the inclusion of the point $0$.  We then have a cartesian square
\begin{equation}\label{eqn:springerbc}
\xymatrix{
\cB \ar[r]^{\iota}\ar[d]_{a} & \tcN \ar[d]^{\mu} \\
\pt \ar[r]_{i_0} & \cN}
\end{equation}

\subsection{Springer correspondence}
\label{subsect:springer}

Let $\bA = \mu_*\bb_\tcN^\flat = \mu_*\bb_\tcN[2d](d)$.  This is a semi\-simple perverse sheaf on $\cN$, known as the \emph{Springer sheaf}.  One approach to studying $\bA$, developed in detail in~\cite{cg:rtcg}, involves \emph{Borel--Moore homology}, which is defined in terms of the hypercohomology of the dualizing complex.  For our purposes, it is convenient to adopt a slightly nonstandard normalization and put
\[
\bmH^i(X) = \bH^i(X, \omega_X[-4d](-2d)).
\]
The Borel--Moore homology of the Steinberg variety $\bmH^\bullet(Z)$ is equipped with a ``convolution product,'' making it into a graded algebra.  Two key results are that there are natural algebra isomorphisms
\begin{equation}\label{eqn:bmh-ext}
\bmH^0(Z) \simeq \Qlb[W]
\qquad\text{and}\qquad
\bmH^\bullet(Z) \simeq \uHom^\bullet(\bA,\bA),
\end{equation}
and that the latter is an isomorphism of graded algebras.  In particular, we have a natural isomorphism $\Qlb[W] \simeq \uEnd(\bA)$, and so an action
\[
\sigma: W \to \Aut(\bA).
\]
Any action of $W$ on $\bA$ would allow us to decompose $\bA$ into isotypic components, but since $\Qlb[W] \simeq \uEnd(\bA)$, we actually have
\begin{equation}\label{eqn:springer-decomp}
\bA \simeq \bigoplus_{\chi \in \Irr(W)} \IC_\chi \otimes V_\chi,
\end{equation}
where the $\IC_\chi$ are various distinct simple perverse sheaves on $\cN$. This labeling of certain simple perverse sheaves by $\Irr(W)$ is what is usually known as the \emph{Springer correspondence}.  For $\chi \in \Irr(W)$, let $\cO_\chi \subset \cN$ be the unique open $G$-orbit in the support of $\IC_\chi$, and let let $L_\chi$ be the irreducible local system given by $\IC_\chi|_{\cO_\chi}[-\dim \cO_\chi](-\frac{1}{2}\dim \cO_{\chi})$.  Thus, $\IC_\chi \simeq \IC(\cO_\chi,L_\chi)$.

\subsection{Modules for $\bmH^\bullet(Z)$}

The convolution product construction also makes the Borel--Moore homology of any subvariety of $\tcN$ into a graded $\bmH^\bullet(Z)$-module.  The convolution action of $\bmH^0(Z)$ on $\bmH^\bullet(\cB) \simeq \bH^\bullet(\cB)^\sharp$ is none other than the action $\kappa$.  

It is clear by base change in~\eqref{eqn:springerbc} that $i_0^*\bA \simeq \bH^\bullet(\cB)^\flat \simeq \bmH^\bullet(\cB)^{\flat\flat}$.  The functor $i_0^*$ therefore induces a map $i_0^*: \uEnd(\bA) \to \End(\bH^\bullet(\cB))$ that is a homomorphism of $\bmH^\bullet(Z)$-modules.  In particular, it is $W$-equivariant, so we have
\[
i_0^*(\sigma(w)f) = \kappa(w) i_0^*(f).
\]

Finally, consider $\bmH^\bullet(\tcN)$.  This is a convolution algebra in its own right.  Since $\bmH^\bullet(\tcN) \simeq \bH^\bullet(\tcN)$, it has another algebra structure given by cup product in ordinary cohomology, but it follows from~\cite[Theorem~8.6.7]{cg:rtcg} (cf.~\cite[Theorem~2.3]{dr:hsv}) that these two algebra structures coincide.  The following theorem of Douglass--R\"ohrle~\cite{dr:hsv} relates the $W$-action on $\bmH^\bullet(\tcN)$ to that on $\bmH^\bullet(Z)$.

\begin{thm}[Douglass--R\"ohrle]\label{thm:dr}
Let $\delta: \tcN \to Z$ be the diagonal embedding.  The induced map $\delta_*: \bmH^\bullet(\tcN) \to \bmH^\bullet(Z)$ in Borel--Moore homology satisfies
\[
\delta_*(w \star f) = w \star \delta_*(f) \star w^{-1}.
\]
for any $w \in W$, where $\star$ denotes the convolution product.
\end{thm}

\subsection{Weakly $\Gm$-equivariant objects}
\label{subsect:wgmeq}

Let $X$ be a variety endowed with an action of $\Gm$.  An object $F \in \Dbc(X)$ is said to be \emph{weakly equivariant} if it is in the image of the forgetful functor $\Forg: \DbcGm(X) \to \Dbc(X)$, where $\DbcGm(X)$ denotes the $\Gm$-equivariant derived category of $X$ in the sense of Bernstein--Lunts~\cite{bl:esf}.  Weakly equivariant objects have the following useful property.

\begin{lem}[{Springer~\cite[Proposition~1]{spr:prfpv}, cf.~Braden \cite[Lemma 6]{bra:hlic}}]\label{lem:lochom}
Let $p: V \to Z$ be a vector bundle, and suppose $\Gm$ acts linearly on the fibers of $p$ with strictly positive weights (or strictly negative weights).  Let $i: Z \to V$ be the inclusion of the zero section.  For a weakly equivariant object $S \in \Dbc(V)$, there are natural isomorphisms $i^!S \to p_!S$ and $p_*S \to i^*S$.\qed
\end{lem}

\begin{rmk}\label{rmk:weq}
Any object obtained by pullback or push-forward of a weakly equivariant object along a $\Gm$-equivariant map is automatically weakly equivariant, and the constant sheaf is always weakly equivariant.  Therefore:
\begin{itemize}
\item If $\Gm$ acts on $\tcN$ by scaling along the the fibers, and on $\cN$ by scaling, then the objects $\bb_{\tcN}$, $\bA$, $\mu^*\bA$, $\mu^!\bA$, and all direct summands of the last three are weakly equivariant.
\item If $\Gm$ acts on $\cN$ by an action that factors through $\tilde G$, where $\tilde G$ is a group isogenous to $G$ with simply-connected derived subgroup, then every semisimple perverse sheaf on $\cN$ is weakly equivariant, since any local system on any nilpotent orbit is $\tilde G$-equivariant.
\end{itemize}
\end{rmk}

\begin{rmk}\label{rmk:allisom}
Suppose $\phi$ is a morphism of functors that is an isomorphism on weakly equivariant objects.  If the domain category of $\phi$ is generated as a triangulated category by weakly equivariant objects, then a standard d\'evissage argument shows that $\phi$ actually induces isomorphisms for \emph{all} objects; in other words, $\phi$ is an isomorphism of functors outright.  This observation will be used when we apply Lemma~\ref{lem:lochom} and other results to $\Dbc(\cN)$, which is generated by the objects of $\PN$.
\end{rmk}

\section{Exactness of $\Phi$}
\label{sect:exact}

In this section, we use hyperbolic localization to prove exactness results for the functor $\Phi = \pi_!\mu^*: \Dbc(\cN) \to \Dbc(\cB)$ of~\eqref{eqn:phi-defn}, as well as for the dual functor
\[
\Phi' = (\pi_*\mu^!)^\sharp: \Dbc(\cN) \to \Dbc(\cB).
\]
To study these functors, we will make use of the additional functor
\[
\Psi = \mu_*\pi^! \simeq (\mu_!\pi^*)^\flat: \Dbc(\cB) \to \Dbc(\cN).
\]
It is clear that $\Psi$ is left-adjoint to $\Phi'$ and right-adjoint to $\Phi$.  In addition, we have
\[
\Phi \simeq \iota^!\mu^*,
\qquad
\Phi' \simeq (\iota^*\mu^!)^\sharp,
\qquad
\Psi(\bb_\cB) \simeq \bA,
\]
with the first two assertions relying on Remark~\ref{rmk:allisom}.  Moreover, since all objects of $\Sh$ are direct sums of copies of $\bb_\cB$, it follows that $\Psi$ restricts to an exact functor of abelian categories
$\Psi: \Sh \to \PN$.  The main result of this section is the following.

\begin{thm}\label{thm:exact}
The functors $\Phi, \Phi': \Dbc(\cN) \to \Dbc(\cB)$ restrict to give isomorphic exact functors of abelian categories $\Phi \simeq \Phi': \PN \to \Sh$.
\end{thm}

We first require the following preliminary result.

\begin{lem}\label{lem:funcisom}
Let $e: \pt \to \cB$ be the inclusion of a point.  There is a natural isomorphism of functors $e^!\pi_*\mu^{!\sharp} \simto e^*\pi_!\mu^*: \Dbc(\cN) \to \Dbc(\pt)$.
\end{lem}
\begin{proof}
For this proof, we fix a choice of Borel subgroup $B \subset G$ and maximal torus $T \subset B$.  Recall that these choices yield a canonical identification $W = N_G(T)/T$.  Let $\fg = \fub \oplus \ft \oplus \fu$  be the corresponding triangular decomposition of $\fg$.  That is, $\ft = \Lie(T)$, $\fu$ is the nilpotent radical of $\Lie(B)$, and $\fub$ is the nilpotent radical of the Lie algebra of the opposite Borel subgroup.

Choose a regular dominant cocharacter $\lambda: \Gm \to T$, and let $\Gm$ act on $\fg$ by composing $\lambda$ with the adjoint action of $T$ on $\fg$.  Clearly, the triangular decomposition of $\fg$ is stable under this action.  Moreover, $\Gm$ acts on $\fu$ with positive weights and on $\fub$ with negative weights, and it acts trivially on $\ft$.  From these observations, it is easy to see that the point $0 \in \cN$ is the unique fixed point for the action of $\Gm$ on $\cN$, and that
\[
\fu = \{ x \in \cN \mid \lim_{\substack{z \in \Gm\\z \to 0}} z \dot x = 0 \}
\qquad\text{and}\qquad
\fub = \{ x \in \cN \mid \lim_{\substack{z \in \Gm\\z \to \infty}} z \dot x = 0 \}.
\]
Consider the following diagram of inclusion maps:
\[
\xymatrix{
\pt \ar[r]^{i}\ar[d]_{\bar\imath} & \fu \ar[d]^{g} \\
\fub \ar[r]_{\bar g} & \cN}
\]
This is a setting in which we may apply the formalism of \emph{hyperbolic localization}, following~\cite{bra:hlic}.  The main theorem of~\cite{bra:hlic} states that there is a natural morphism of functors $\bar\imath^* \bar g^! \to i^! g^*$ that is an isomorphism on weakly equivariant objects.  By Remark~\ref{rmk:allisom}, this is an isomorphism of functors in our situation.

Next, consider the constant maps $p: \fu \to \pt$ and $\bar p: \fub \to \pt$.  Using Lemma~\ref{lem:lochom}, we obtain a natural isomorphism $\bar p_*\bar g^! \overset{\sim}{\to} p_!g^*$.  Finally, let $e: \pt \to \cB$ (resp.~$\bar e: \pt \to \cB$) be the inclusion of the point corresponding to the Borel subgroup $B$ (resp.~the opposite Borel subgroup to $B$).  Forming pullbacks over $\pi$, we obtain the diagrams
\[
\xymatrix{
 & \fu \ar[d]_{\tilde e}\ar[r]^{p} &\pt \ar[d]^e \\
\cN & \tcN \ar[l]_{\mu} \ar[r]^{\pi} & \cB}
\qquad
\xymatrix{
 & \fub \ar[d]_{\Tilde{\Bar{e}}}\ar[r]^{\bar p} &\pt \ar[d]^{\bar e} \\
\cN & \tcN \ar[l]_{\mu} \ar[r]^{\pi} & \cB}
\]
It is clear that $\mu\tilde e = g$ and $\mu\Tilde{\Bar{e}} = \bar g$.  By base change, we have $p_!\tilde e^* \simeq e^*\pi_!$ and $\bar p_*{\Tilde{\Bar{e}}}^! \simeq {\bar e}^!\pi_*$.  Combining with our earlier observations, we obtain an isomorphism $e^!\pi_*\mu^! \simto \bar e^* \pi_! \mu^*$.  It is clear that on $\Dbc(\cB)$, we have $\bar e^* \simeq e^*$ and $e^! \simeq e^{*\sharp}$, so we now have the isomorphism
$e^*\pi_*\mu^{!\sharp} \simto e^* \pi_!\mu^*$, as desired.
\end{proof}

\begin{proof}[Proof of Theorem~\ref{thm:exact}]
We begin by showing that $\Phi|_{\PN}$ is exact.  Let $F \in \PN$.  We wish to show that $H^j(\Phi(F)) = 0$ for $j \ne 0$.  First, observe that for $n < 0$, we have
\[
\uHom(\Phi(F), \bb_\cB[n]) \simeq \uHom(F, \Psi(\bb_\cB)[n]) \simeq \uHom(F, \bA[n]) = 0.
\]
This shows that $H^j(\Phi(F)) = 0$ for $j > 0$.

By duality, we have that $H^j(\pi_*\mu^!F) \simeq H^j(\D \pi_!\mu^*\D F) \simeq (\D H^{-2d-j}(\pi_!\mu^*\D F))^\flat$.  Thus, $H^j(\pi_*\mu^!F) = 0$ if $j < -2d$.  For the inclusion of a point $e: \pt \to \cB$, we know that $e^*$ is an exact functor on $\Sh$, so it follows that $H^j(e^*\pi_*\mu^!F^\sharp) = 0$ if $j < 0$.  By Lemma~\ref{lem:funcisom}, this implies that $H^j(e^*\pi_!\mu^*F) = 0$ for $j < 0$ as well.  Now, $e^*$ is also faithful on the abelian category $\Sh$ (though not, of course, on $\Dbc(\cB)$), so that vanishing implies that $H^j(\pi_!\mu^*F) = 0$ for $j < 0$, as desired.  Thus, $\Phi|_{\PN}$ is exact.

Finally, note that the exact functor $e^*: \Sh \to \Shpt$ is an equivalence of categories.  The isomorphism $\Phi|_{\PN} \simeq \Phi'|_{\PN}$ then follows from Lemma~\ref{lem:funcisom}.
\end{proof}

\begin{cor}\label{cor:pure}
The functor $\Phi \simeq \Phi': \PN \to \Sh$ preserves purity.  That is, it takes pure objects of weight $w$ to pure objects of weight $w$.
\end{cor}
\begin{proof}
Let $\cF \in \PN$ be pure of weight $w$.  By the well-known rules~\cite[5.1.14]{bbd} for behavior of weights under various sheaf functors, we see that $\Phi = \pi_!\mu^*$ takes $\cF$ to an object with weights${}\le w$, whereas $\Phi' = \pi_*\mu^!$ takes it to one with weights${}\ge w$.
\end{proof}

\section{Action of the Weyl Group}
\label{sect:weyl}

In Section~\ref{sect:prelim}, we considered the action $\sigma$ of $W$ on $\uEnd(\bA)$, and the action $\kappa$ on $\bH^\bullet(\cB)$.  In this section, we discuss several additional actions, and prove a $W$-equivariance result for $\Psi$.  There are two natural commuting actions $\lambda,\rho: W \to \Aut(\uHom^i(\bA,\bA))$, given by
\[
\lambda(w)(f) = \sigma(w) \circ f
\qquad\text{and}\qquad
\rho(w)(f) = f \circ \sigma(w^{-1}).
\]
The exactness result of Section~\ref{sect:exact} allows us to construct a new action as follows.

\begin{prop}\label{prop:W-action}
For any $F \in \PN$, the sheaf $\Phi(F)$ carries a natural action of $W$.  If $F$ is simple, then we have
\[
\Phi(F) \simeq
\begin{cases}
0 & \text{if $F \notin \Spr$,} \\
\bb_{\cB} \otimes V_\chi^* & \text{if $F \simeq \IC_\chi$.}
\end{cases}
\]
\end{prop}
\begin{proof}
The following general principle is easy to see: if $\cA$ is a semisimple $\Bbbk$-linear abelian category containing a unique simple object $S$ up to isomorphism, and $\End(S) \simeq \Bbbk$, then any object $A$ is canonically isomorphic to $S \otimes \Hom(S,A)$.  Applying this to $\Sh$, we have $\Phi(F) \simeq \bb_{\cB} \otimes \uHom(\bb_{\cB}, \Phi(F))$.  By adjunction, $\uHom(\bb_{\cB}, \Phi(F)) \simeq \uHom(\Psi(\bb_{\cB}),F) \simeq \uHom(\bA,F)$. 

The $W$-action on $\bA$ induces one on $\uHom(\bA,F)$ for any $F$, and therefore on $\Phi(F) \simeq \bb_{\cB} \otimes \uHom(\bA,F)$.  For simple $F$, it is clear from~\eqref{eqn:springer-decomp} that $\Hom(\bA,F) = 0$ if $F \notin \Spr$, and that $\Hom(\bA,\IC_\chi) \simeq V_\chi^*$.
\end{proof}

The action described in this proposition will be denoted $\nu: W \to \Aut(\Phi(F))$.  This action gives rise a $W$-action on the vector space $\Gamma(\Phi(F))$ for any $F \in \PN$.  In other words, the functor $\Gamma \circ \Phi$ may be regarded as taking values in $\Rep(W)$.  Since $\PN$ and $\Rep(W)$ are both semisimple categories, the following result is immediate from Proposition~\ref{prop:W-action}.

\begin{thm}
For $F \in \PN$, the functor $\Gamma\circ\Phi: \PN \to \Rep(W)$ is given by
\[
(\Gamma \circ \Phi)(F) \simeq
\begin{cases}
0 & \text{if $F \notin \Spr$,} \\
V_\chi^* & \text{if $F \simeq \IC_\chi$.}
\end{cases}
\]
In particular, $\Gamma\circ\Phi|_{\Spr}: \Spr \to \Rep(W)$ is an equivalence of categories.\qed
\end{thm}

\begin{cor}\label{cor:duality}
The category $\Spr \subset \PN$ is stable under Verdier duality $\D$.  In fact, for each simple perverse sheaf $\IC_\chi \in \Spr$, we have $\D\IC_\chi \simeq \IC_\chi$.
\end{cor}
\begin{proof}
It follows from Theorem~\ref{thm:exact} and the fact that $\cB$ is a projective variety that $\Gamma \circ \Phi$ commutes with $\D$.  For a simple perverse sheaf $F \in \PN$, we see that $(\Gamma\circ\Phi)(F) \ne 0$ if and only if $(\Gamma \circ\Phi)(\D F) \ne 0$, so $\D$ preserves $\Spr$.  Moreover, for $F \simeq \IC_\chi$, we have $(\Gamma \circ \Phi)(\D\IC_\chi) \simeq \D(V_\chi^*) \simeq V_\chi$.  The result follows using the noncanonical isomorphism~\eqref{eqn:contragr}.
\end{proof}

Note that when $F = \bA$, $\nu$ is obtained via the adjunction isomorphism
\begin{equation}\label{eqn:end-adj}
\theta: \uHom(\bb_\cB,\Phi(\bA)) \simto \uHom(\bA,\bA)
\end{equation}
from the action on $\Hom(\bA,\bA)$ that we have called $\rho$.  The other action $\lambda$ on $\Hom(\bA,\bA)$ also induces an action on $\Phi(\bA)$, which we denote
\[
\hat\lambda: W \to \Aut(\Phi(\bA)).
\]
Since $\rho$ and $\lambda$ commute, $\nu$ and $\hat\lambda$ commute as well.
By an abuse of notation, we will also write $\nu$ and $\hat\lambda$ for the corresponding actions on the vector space $\uHom^i(\bb_\cB, \Phi(\bA))$.  By definition, we have
\begin{equation}\label{eqn:theta-eq0}
\theta(\hat\lambda(v)\nu(w)f) = \lambda(v)\rho(w) \theta(f) = \sigma(v) \circ \theta(f) \circ \sigma(w^{-1}).
\end{equation}

\begin{lem}\label{lem:dr-equiv}
The map $\Psi: \uHom^i(\bb_\cB,\bb_\cB) \to \uHom^i(\bA,\bA)$ has the property that $\Psi(\kappa(w)f) = \sigma(w) \circ \Psi(f) \circ \sigma(w^{-1})$.
\end{lem}
\begin{proof}
This is essentially a restatement of Theorem~\ref{thm:dr} due to Douglass--R\"ohrle.  Note first that the functor $\pi^!$ induces a $W$-equivariant isomorphism $\bH^\bullet(\cB) \simto \bH^\bullet(\tcN)$, so it suffices to study the $W$-equivariance of $\mu_*: \Dbc(\tcN) \to \Dbc(\cN)$.  Let $q: Z \to \cN$ be the natural projection map.  We then have a commutative diagram:
\[
\xymatrix{
& \tcN \ar[dl]_\delta \ar[dr]^\mu \\
Z \ar[rr]_q && \cN }
\]
Recall~\cite[Section 8.6]{cg:rtcg} that the second isomorphism in~\eqref{eqn:bmh-ext} arises from the fact that $q_*\omega_Z^{\sharp\sharp} \simeq \cRHom(\bA,\bA)$.  There is also a natural adjunction map $\delta_*\omega_\tcN \to \omega_Z$.  Consider the composition
\[
\xymatrix{
\mu_*\omega_\tcN^{\sharp\sharp} \ar[r] & q_*\omega_Z^{\sharp\sharp} \ar[r]^-\sim
& \cRHom(\bA,\bA).}
\]
Applying $R\Gamma$ to the first map yields the induced map $\delta_*: \bmH^\bullet(\tcN) \to \bmH^\bullet(Z)$ in Borel--Moore homology.  On the other hand, we can identify $\omega_\tcN^{\sharp\sharp} \simeq \bb_\tcN \simeq \cRHom(\bb_\tcN,\bb_\tcN)$, and then the composition above becomes the canonical morphism
\[
\mu_*\cRHom(\bb_\tcN,\bb_\tcN) \to \cRHom(\bA,\bA),
\]
and applying $R\Gamma$ gives us the map $\uHom^\bullet(\bb_\tcN,\bb_\tcN) \to \uHom^\bullet(\bA,\bA)$ induced by $\mu$.  To summarize, we have the following commutative diagram:
\[
\xymatrix{
\bmH^\bullet(\tcN) \ar[r]^-\sim \ar[d]_{\delta_*} &
\uHom^\bullet(\tcN,\tcN) \ar[d]^{\mu} \\
\bmH^\bullet(Z) \ar[r]_-q &
\uHom^\bullet(\bA,\bA) }
\]
Since the top isomorphism is $W$-equivariant, and the bottom one is an algebra isomorphism sending $w \in \Qlb[W] \simeq \bmH^0(Z)$ to $\sigma(w) \in \uEnd(\bA)$, the result follows from Theorem~\ref{thm:dr}.
\end{proof}

\begin{prop}\label{prop:ext-tensor}
There is a natural isomorphism $\alpha: \bH^\bullet(\cB) \otimes \uHom(\bb_\cB,\Phi(\bA)) \simto \uHom^\bullet(\bb_\cB,\Phi(\bA))$.  Its composition with the adjunction $\theta$, denoted
\[
\Theta = \theta\circ\alpha: \bH^\bullet(\cB) \otimes \uHom(\bb_\cB,\Phi(\bA)) \simto \uHom^\bullet(\bA,\bA),
\]
is $W$-equivariant in the following way: for $u,v,w \in W$, we have
\[
\Theta( \kappa(u)f \otimes \hat\lambda(v)\nu(w)g) = 
\lambda(v)\rho(w)\Theta(\kappa(w^{-1}u)f \otimes g).
\]
\end{prop}
\begin{proof}
Recall from the proof of Proposition~\ref{prop:W-action} that $\Phi(\bA) \simeq \bb_\cB \otimes \uHom(\bb_\cB,\Phi(\bA))$.  It follows that 
\begin{multline*}
\uHom^i(\bb_\cB, \Phi(\bA)) \simeq
\uHom^i(\bb_\cB, \bb_\cB \otimes \uHom(\bb_\cB,\Phi(\bA))) \\
\simeq \uHom^i(\bb_\cB, \bb_\cB) \otimes \uHom(\bb_\cB,\Phi(\bA)).
\end{multline*}
Since $\bH^i(\cB) \simeq \uHom^i(\bb_\cB,\bb_\cB)$, we obtain the desired isomorphism $\alpha$.  Note that $\alpha$ is given by composition: that is, if $f \in \uHom^i(\bb_\cB,\bb_\cB)$ and $g \in \uHom(\bb_\cB,\Phi(\bA))$, then
\[
\alpha(f \otimes g) = g \circ f.
\]
The adjunction isomorphism $\theta$ behaves on compositions according to the rule $\theta(g \circ f) = \theta(g) \circ \Psi(f)$.  Using~\eqref{eqn:theta-eq0} and Lemma~\ref{lem:dr-equiv}, we find
\begin{align*}
\Theta(\kappa(u)f \otimes \hat\lambda(v)\nu(w)g)
&= \theta(\hat\lambda(v)\nu(w)g) \circ \Psi(\kappa(w)f) \\
&= \sigma(v) \circ \theta(g) \sigma(w^{-1}) \circ \sigma(u) \circ \Psi(f) \circ \sigma(u^{-1}) \\
&= \sigma(v) \circ \theta(g) \circ \sigma(w^{-1}u) \circ \Psi(f) \circ \sigma(u^{-1}w) \circ \sigma(w^{-1}) \\
&= \sigma(v) \circ \theta(g) \circ \Psi(\kappa(w^{-1}u)f) \circ \sigma(w^{-1}) \\
&= \lambda(v)\rho(w) \theta(g \circ \kappa(w^{-1}u)f).
\end{align*}
The result follows.
\end{proof}

\begin{thm}\label{thm:ext-w-invt}
For $\chi,\psi \in \Irr(W)$, there is a natural isomorphism
\begin{equation}\label{eqn:ext-phi}
\uHom^i(\Phi(\IC_\chi), \Phi(\IC_\psi)) \simeq V_\chi \otimes \bH^i(\cB) \otimes V_\psi^*,
\end{equation}
and thus $\uHom^i(\Phi(\IC_\chi), \Phi(\IC_\psi))$ is endowed with a natural action of $W$.  Moreover, $\Phi$ induces isomorphisms
\[
\uHom^i(\IC_\chi,\IC_\psi) \simeq
\uHom^i(\Phi(\IC_\chi), \Phi(\IC_\psi))^W \simeq (V_\chi \otimes \bH^i(\cB) \otimes V_\psi^*)^W.
\]
\end{thm}
\begin{proof}
The isomorphism~\eqref{eqn:ext-phi} is immediate from Proposition~\ref{prop:W-action}.  Next, using~\eqref{eqn:springer-decomp}, we can decompose $\uHom^i(\bA,\bA)$ as
\[
\Hom^i(\bA,\bA)
\simeq \bigoplus_{\chi,\psi} \Hom^i(\IC_\chi,\IC_\psi) \otimes V_\chi^* \otimes V_\psi.
\]
Thus, in terms of the action of $W \times W$ on $\Hom^i(\bA,\bA)$ by $\lambda \boxtimes \rho$, we can find $\Hom^i(\IC_\chi,\IC_\psi)$ by picking out the $\chi^* \boxtimes \psi$-isotypic component:
\[
\Hom^i(\IC_\chi,\IC_\psi) \simeq \Hom_{W\times W}(V_\chi^* \boxtimes V_\psi, \Hom^i(\bA,\bA)).
\]
Using Proposition~\ref{prop:ext-tensor}, this is isomorphic to
\[
\uHom_{W \times W}(V_\chi^* \boxtimes V_\psi, \bH^i(\bb_{\cB}) \otimes \uHom(\bb_{\cB},\Phi(\bA))),
\]
where $W \times W$ acts on $\bH^i(\bb_{\cB}) \otimes \uHom(\bb_{\cB},\Phi(\bA))$ by $\hat\lambda \boxtimes (\kappa \otimes \nu)$.  That is, for $v,w \in W$ and $f \otimes g \in \bH^i(\bb_{\cB}) \otimes \uHom(\bb_{\cB},\Phi(\bA))$, we put
\[
(v,w)\cdot(f \otimes g) = \kappa(w)f \otimes \hat\lambda(v)\nu(w)g.
\]
Using the adjunction~\eqref{eqn:end-adj} and the isomorphism~\eqref{eqn:bmh-ext}, we see that $\uHom(\bb_{\cB},\Phi(\bA))$ decomposes under $\hat\lambda \boxtimes \nu$ as $\uHom(\bb_\cB, \Phi(\bA)) \simeq \bigoplus_\phi V_\phi \otimes V_\phi^*$.  Picking off the $\chi^*$-isotypic component for the first factor of $W$, we find that
\[
\uHom^i(\IC_\chi,\IC_\psi) \simeq \Hom_W(V_\psi, \bH^i(\cB) \otimes  V_\chi))
\simeq (V_\psi^* \otimes \bH^i(\cB) \otimes V_\chi)^W.\qedhere
\]
\end{proof}

\begin{cor}\label{cor:ext-purity}
For $\chi,\psi \in \Irr(W)$, $\uHom^i(\IC_\chi,\IC_\psi)$ vanishes if $i$ is odd, and is pure of weight $i$ if $i$ is even.
\end{cor}
\begin{proof}
This follows from the previous theorem and the well-known fact that $\bH^i(\cB)$ vanishes if $i$ is odd and is pure of weight $i$ if $i$ is even.
\end{proof}

\section{Orthogonal Decomposition of $\Dbc(\cN)$}
\label{sect:orth}

For a $G$-stable locally closed subvariety $Y \subset \cN$, and let $\Dspr(Y) \subset \Dbc(Y)$ be the full triangulated subcategory generated by objects $\IC_\chi|_Y$ with $\cO_\chi \subset Y$.  On the other hand, let $\Dsprp(Y) \subset \Dbc(Y)$ be the full triangulated subcategory generated by simple perverse sheaves $\IC(\cO,L)|_Y$ with $\cO \subset Y$ but $\IC(\cO,L) \notin \Spr$.

\begin{thm}
For any $G$-stable locally closed subvariety $u: Y \hookrightarrow \cN$, we have
\begin{equation}\label{eqn:orth}
\Dbc(Y) \simeq \Dspr(Y) \oplus \Dsprp(Y).
\end{equation}
Moreover, if $s: Z \hookrightarrow Y$ is the inclusion of a smaller $G$-stable locally closed subvariety, the functors $s^*$ and $s^!$ respect this decomposition: we have
\begin{equation}\label{eqn:orth-stab}
\begin{aligned}
s^*(\Dspr(Y)),\ s^!(\Dspr(Y)) &\subset \Dspr(Z), \\
s^*(\Dsprp(Y)),\ s^!(\Dsprp(Y)) &\subset \Dsprp(Z).
\end{aligned}
\end{equation}
\end{thm}
\begin{proof}
If $s: Z \to Y$ is an open embedding, then~\eqref{eqn:orth-stab} is obvious.  Since the inclusion of any locally closed subvariety can be factored as a closed embedding followed by an open embedding, we henceforth treat~\eqref{eqn:orth-stab} only in the case where $Z$ is closed in $Y$.

Let $n_Y$ denote the number of nilpotent orbits in $\overline{Y} \smallsetminus Y$.  We will prove~\eqref{eqn:orth} and~\eqref{eqn:orth-stab} simultaneously by induction on $n_Y$.  Note that~\eqref{eqn:orth} is equivalent to the assertion that for $F \in \Dspr(Y)$ and $F' \in \Dsprp(Y)$, we have $\uHom(F,F') = \uHom(F',F) = 0$.  We can further reduce to the case where $F$ and $F'$ are shifts of simple perverse sheaves.  That is,~\eqref{eqn:orth} is equivalent to the statement that if $F$ and $F'$ are simple perverse sheaves with $F \in \Dspr(Y)$ and $F' \in \Dsprp(Y)$, then
\begin{equation}\label{eqn:orth-van}
\uHom^i(F,F') = \uHom^i(F',F) = 0
\qquad\text{for all $i \ge 0$.}
\end{equation}

We begin by proving~\eqref{eqn:orth-van} in the case where $n_Y = 0$, i.e., when $Y$ is closed in $\cN$.  In fact, since $u_*: \Dbc(Y) \to \Dbc(\cN)$ is faithful for any closed $Y \subset \cN$, we may reduce to the case where $Y = \cN$.  Since $F \in \Spr$, $F$ is a direct summand of $\bA$, and it suffices to show that $\uHom^i(\bA,F') = \uHom^i(F',\bA) = 0$.  Since $\bA \simeq \Psi(\bb_\cB)$, we have by adjunction that
\[
\uHom^i(\bA,F') \simeq \uHom^i(\bb_\cB, \Phi'(F'))
\qquad\text{and}\qquad
\uHom^i(F',\bA) \simeq \uHom^i(\Phi(F),\bb_\cB).
\]
Since $\Phi(F') = \Phi'(F') = 0$ by Theorem~\ref{thm:exact} and Proposition~\ref{prop:W-action}, the desired vanishing holds.

Suppose now that~\eqref{eqn:orth} is known to hold whenever $n_Y \le k$.  Let us prove~\eqref{eqn:orth-stab}.  Since $Z$ is assumed to be a closed subvariety of $Y$, we clearly have $n_Z \le n_Y$; in particular, we know that $\Dbc(Z) \simeq \Dspr(Z) \oplus \Dsprp(Z)$.  Therefore, given $F \in \Dspr(Y)$, we have a canonical decomposition $s^*F \simeq (s^*F)_{\Spr} \oplus (s^*F)_{\Spr}^\perp$ with $(s^*F)_{\Spr} \in \Dspr(Z)$ and $(s^*F)_{\Spr}^\perp \in \Dsprp(Z)$.  We wish to prove that $(s^*F)_{\Spr}^\perp = 0$.  If that is not the case, there certainly exists some simple perverse sheaf $F' \in \Dsprp(Z)$ and some $i \in \Z$ such that $\uHom^i((s^*F)_{\Spr}^\perp, F') \ne 0$.  We also know that $\uHom^\bullet((s^*F)_{\Spr}, F') = 0$, so it follows that
\[
0 \ne \uHom^i((s^*F)_{\Spr}^\perp, F') \simeq \uHom^i(s^*F,F') \simeq \uHom^i(F,s_*F').
\]
But $s_*F'$ is clearly a simple perverse sheaf in $\Dsprp(Y)$, and since $F \in \Dspr(Y)$, we have a contradiction.  Thus, $s^*F \in \Dspr(Z)$.  The proofs of the other assertions in~\eqref{eqn:orth-stab} are parallel.

Now, suppose that~\eqref{eqn:orth} and~\eqref{eqn:orth-stab} are both known for $n_Y \le k$, and let us prove~\eqref{eqn:orth-van} when $n_Y = k+1$.  Let $F$ and $F'$ be simple perverse sheaves on $Y$ with $F \in \Dspr(Y)$ and $F' \in \Dsprp(Y)$.  Choose an orbit $\cO$ that is open in $\overline{Y} \smallsetminus Y$, and let $\tilde Y = Y \cup \cO$.  Let $s: \cO \to \tilde Y$ and $j: Y \to \tilde Y$ be the inclusion maps, and consider the objects $j_{!*}F, j_{!*}F' \in \Dbc(\tilde Y)$.  These are simple perverse sheaves on $\tilde Y$; moreover, we clearly have $j_{!*}F \in \Dspr(\tilde Y)$ and $j_{!*}F' \in \Dsprp(\tilde Y)$.  Form the distinguished triangle
\[
s_*s^! j_{!*}F' \to j_{!*}F' \to j_*j^*j_{!*}F' \to.
\]
Note that $j^*j_{!*}F' \simeq F'$.  Next, form the long exact sequence
\[
\cdots \to \uHom^i(j_{!*}F, j_{!*}F') \to \uHom^i(j_{!*}F, j_*F') \to \uHom^{i+1}(j_{!*}F',s_*s^!j_{!*}F') \to \cdots.
\]
Note that $n_{\tilde Y} = n_Y - 1$, so~\eqref{eqn:orth-van} holds on $\tilde Y$ by assumption.  In particular, we have $\uHom^\bullet(j_{!*}F, j_{!*}F') = 0$.  Since~\eqref{eqn:orth-stab} also holds by assumption, we have $s^*j_{!*}F \in \Dspr(\cO)$ and $s^!j_{!*}F' \in \Dsprp(\cO)$, so
\[
\uHom^\bullet(j_{!*}F',s_*s^!j_{!*}F') =
\uHom^\bullet(s^*j_{!*}F', s^!j_{!*}F') = 0.
\]
We conclude that $\uHom^i(j_{!*}F,j_*F') = 0$.  But that means $\uHom^i(j^*j_{!*}F,F') \simeq \uHom^i(F,F') = 0$ as well, as desired.
\end{proof}

One concrete consequence of this theorem is the following.

\begin{cor}\label{cor:spr-locsys}
If $L$ is an irreducible local system on an orbit $\cO \subset \cN$ that occurs as a composition factor of some cohomology sheaf $H^i(\IC_\chi|_{\cO})$ with $\chi \in \Irr(W)$, then $L \simeq L_\psi$ for some $\psi \in \Irr(W)$. \qed
\end{cor}

\section{Green Functions}
\label{sect:green}

The aim of this section is to study the restrictions $\IC_\chi|_{\cO}$ as $\cO$ varies over the $G$-orbits in $\cN$.  Specifically, we encode information about these restrictions in a family of polynomials $p_{\chi,\psi}(t)$, sometimes called \emph{Green functions}.  The main result, Theorem~\ref{thm:ls}, gives a way to compute these polynomials, by relating them to the known groups $\Hom^i(\IC_\chi,\IC_\psi)$.

We begin with some notation.  For a variety $X$, let $K(X)$ denote the quotient of the Grothendieck group of $\Dbc(X)$ obtained by identifying the classes of simple perverse sheaves of the same weight that become isomorphic when the Weil structure is forgotten.  (Thus, $K(X)$ does not detect twists in the Weil structure by a root of unity.)  For an object $F \in \Dbc(X)$, let $[F]$ denote its class in $K(\cN)$.  This Grothendieck group is naturally a module over the Laurent polynomial ring $\Z[t,t^{-1}]$, where multiplication by $t$ corresponds to Tate twist: $[F(1)] = t^{-1}[F]$.  By choosing a square root of the Tate sheaf, we can regard $K(X)$ as a module over $\Z[t^{1/2},t^{-1/2}]$, with $[F(\frac{1}{2})] = t^{-1/2}[F]$.  For instance, the group $K(\pt)$ is a free $\Z[t^{1/2},t^{-1/2}]$-module of rank $1$, generated by the class $[\bb_\pt]$ of a $1$-dimensional vector space of weight $0$.  

For an orbit $\cO \subset \cN$, the group $K(\cO)$ is a free $\Z[t^{1/2},t^{-1/2}]$-module generated by the classes of irreducible local systems on $\cO$.  In view of Corollary~\ref{cor:spr-locsys}, we may write
\begin{equation}\label{eqn:p-defn}
[\IC_\chi|_{\cO}] = \sum_{\{\psi \mid \cO_\psi = \cO\}} p_{\chi,\psi}(t)[L_\psi]
\qquad\text{for some $p_{\chi,\psi}(t) \in \Z[t^{1/2},t^{-1/2}]$}.
\end{equation}
(That is, only local systems belonging to $\Dspr(\cO)$ may appear on the right-hand side.)  We clearly have $\IC_\chi|_{\cO_\chi} \simeq L_\chi[\dim \cO_\chi](\frac{1}{2}\dim \cO_\chi)$ and $\IC_\chi|_{\cO} = 0$ if $\cO \not\subset \overline{\cO_\chi}$.  In other words:
\begin{equation}\label{eqn:p-cond}
p_{\chi,\psi}(t) =
\begin{cases}
t^{-(\dim \cO_\chi)/2} & \text{if $\psi = \chi$,} \\
0 & \text{if $\cO_\psi \not\subset \overline{\cO_\chi}$, or if $\cO_\psi = \cO_\chi$ and $\psi \ne \chi$.}
\end{cases}
\end{equation}
Our goal is to determine the polynomials $p_{\chi,\psi}(t)$.  

Another description of these polynomials is as follows.  Each cohomology sheaf $H^i(\IC_\chi|_{\cO})$ is, of course, a finite-length object in the category of local systems on $\cO$.  Let $(H^i(\IC_\chi|_{\cO}) : L_\psi(j))$ denote the multiplicity of $L_\psi(j)$ in any composition series.  We then have
\begin{equation}\label{eqn:p-cohom}
p_{\chi,\psi}(t) = \sum_{i \in \Z,\ j \in \frac12\Z} (-1)^i (H^i(\IC_\chi|_{\cO}) : L_\psi(j)) t^{-j}.
\end{equation}
A result of Springer leads to a tremendous simplification of this formula; see Remark~\ref{rmk:spr-pure}.  Actually, $p_{\chi,\psi}(t)$ lies in $\Z[t^{-1}]$ and has nonnegative coefficients (see~\eqref{eqn:p-pure}), but we will not require these facts.

To state the main result, we require two additional families of polynomials. First, define $\lambda_{\chi,\psi}(t) \in \Z[t^{1/2},t^{-1/2}]$ as follows:
\begin{align}
[R\Gamma_c(\cO_\chi, L_\chi \otimes L_\psi)] &= \lambda_{\chi,\psi}(t)[\bb_{\pt}] &&\text{if $\cO_\chi = \cO_\psi$,} \notag\\
\lambda_{\chi,\psi}(t) &= 0 &&\text{if $\cO_\chi \ne \cO_\psi$.} \label{eqn:lam-cond}
\end{align}
Second, define polynomials $\omega_{\chi,\psi}(t)$ by
\[
[\D \RHom(\IC_\chi,\IC_\psi)] = \omega_{\chi,\psi}(t)[\bb_{\pt}].
\]
Recall from Corollary~\ref{cor:ext-purity} that $\RHom(\IC_\chi,\IC_\psi)$ is pure of weight $0$ and has vanishing cohomology in odd degrees.  The same statements then hold for its dual as well.  As with $p_{\chi,\psi}(t)$ in~\eqref{eqn:p-pure}, it follows that
\[
\omega_{\chi,\psi}(t) = \sum_{i \in \Z} \dim H^{2i}(\D\RHom(\IC_\chi,\IC_\psi))t^i = \sum_{i \in \Z} \dim \uHom^{-2i}(\IC_\chi,\IC_\psi) t^i.
\]
The coinvariant algebra of $W$ has the property that the $W$-action in complementary degrees is related by tensoring with the sign character $\varepsilon$: that is, $\bH^i(\cB) \simeq \bH^{2d-i}(\cB) \otimes \varepsilon$ as $W$-representations.  Using the noncanonical isomorphism~\eqref{eqn:contragr} together with Theorem~\ref{thm:ext-w-invt}, we can rewrite the above formula as
\begin{equation}\label{eqn:omega-coinvt}
\omega_{\chi,\psi}(t) = t^{-2d}\sum_{i \in \Z} \dim \Hom_W(V_\chi \otimes V_\psi \otimes \varepsilon, \bH^{2i}(\cB)) t^i.
\end{equation}
The main result of this section is the following.

\begin{thm}\label{thm:ls}
The matrices $P = (p_{\chi,\psi})$, $\Lambda = (\lambda_{\chi,\psi})$, and $\Omega = (\omega_{\chi,\psi})$ satisfy 
\begin{equation}\label{eqn:ls}
P \Lambda P^{\mathrm t} = \Omega,
\end{equation}
where $P^{\mathrm t}$ is the transpose of $P$.
In other words, given $\chi,\psi \in \Irr(W)$, we have
\[
\sum_{\phi,\phi' \in \Irr(W)} p_{\chi,\phi}(t)\lambda_{\phi,\phi'}(t)p_{\psi,\phi'}(t) = \omega_{\chi,\psi}(t).
\]
Moreover, $P$ and $\Lambda$ are the unique matrices with entries in $\Q(t^{1/2})$ satisfying~\eqref{eqn:ls}, \eqref{eqn:p-cond}, and~\eqref{eqn:lam-cond}.
\end{thm}

\begin{rmk}\label{rmk:spr-pure}
This theorem is essentially equivalent to the part of~\cite[Theorem~24.8]{lus:cs5} relevant to $\Spr$.  The most substantial difference is that in {\it loc.~cit.}, the polynomials $p_{\chi,\psi}(t)$ are defined in a slightly different way.  Correcting for different normalization conventions (see Remark~\ref{rmk:ls-alg}), the definition in~\cite{lus:cs5} is
\begin{equation}\label{eqn:p-pure}
p_{\chi,\psi}(t) = \sum_{i \in \Z} (H^{2i}(\IC_\chi|_{\cO}) : L_\psi(-i)) t^i.
\end{equation}
The equivalence of this formula with~\eqref{eqn:p-cohom} is implied by an important result of Springer~\cite{spr:prfpv}, which states that for any $\chi \in \Irr(W)$ and any orbit $\cO \subset \cN$, the object $\IC_\chi|_{\cO} \in \Dbc(\cO)$ is pure of weight $0$, and that $H^i(\IC_\chi|_{\cO}) = 0$ if $i$ is odd.  The proof of~\eqref{eqn:ls} in~\cite{lus:cs5} also relies on Springer's purity theorem.
\end{rmk}

\begin{rmk}[Lusztig--Shoji algorithm]\label{rmk:ls-alg}
The uniqueness asserted in Theorem~\ref{thm:ls} is proved by Lusztig~\cite{lus:cs5} in a very explicit constructive way.  This proof, which will not be repeated here, consists primarily of a description of an algorithm for finding $P$ and $\Lambda$ from knowledge of $\Omega$.  Since $\Omega$ can be described as in~\eqref{eqn:omega-coinvt} using only the representation theory of $W$, this algorithm can be effectively used to compute the $p_{\chi,\psi}(t)$.  Generalizations of this algorithm, sometimes called the \emph{Lusztig--Shoji algorithm}, have been studied in~\cite{aa:scdg, ah:ocenc, gm:sp, sho:gf1, sho:gf2, sho:gfls}, and a computer implementation is available at~\cite{a:iglsa}.

The reader should be aware that Lusztig originally used polynomials $p'_{\chi,\psi}(t)$ and $\omega'_{\chi,\psi}(t)$ following different normalization conventions, while the recent works mentioned above involve polynomials $p''_{\chi,\psi}(t)$ and $\omega''_{\chi,\psi}(t)$ following a third convention.  The relationship among these is as follows:
\begin{align*}
p'_{\chi,\psi}(t) &= t^{\frac{1}{2}\dim \cO_\chi}p_{\chi,\psi}(t) &
p''_{\chi,\psi}(t) &= t^d p_{\chi,\psi}(t) \\
\omega'_{\chi,\psi}(t) &= t^{\frac{1}{2}(\dim \cO_\chi + \dim \cO_\psi)} \omega_{\chi,\psi}(t) &
\omega''_{\chi,\psi}(t) &= t^{2d} \omega_{\chi,\psi}(t).
\end{align*}
\end{rmk}

For the next three lemmas, let $j_{\cO}: \cO \to \cN$ denote the inclusion of an orbit.

\begin{lem}\label{lem:orbit-sum}
For any $F, F' \in \Dbc(\cN)$, we have
\[
[\D\RHom(F,F')] = \sum_{\cO \subset \cN} [\D\RHom(j_{\cO}^*F, j_{\cO}^!F')].
\]
\end{lem}
\begin{proof}
Let $X \subset \cN$ denote the closure of the support of $F'$, and let $\cO_0 \subset X$ be an orbit that is open in $X$.  Let $h: Y \to \cN$ be the inclusion of the closed subset $Y = X \smallsetminus \cO_0$, and let $s: (\cN \smallsetminus Y) \to \cN$ be the inclusion of its open complement. Consider the distinguished triangle $h_*h^!F' \to F' \to s_*s^*F' \to$.  Since $s^*F' \simeq s^!F'$ is supported on $\cO$, we see that $s_*s^*F'$ is naturally isomorphic to $j_{\cO_0*}j_{\cO_0}^!F'$.  Applying $\D\RHom(F,\cdot)$ and the usual adjunction properties, we obtain a distinguished triangle
\[
\D\RHom(j_{\cO_0}^*F, j_{\cO_0}^!F') \to \D\RHom(F,F') \to \D\RHom(h^*F, h^!F') \to,
\]
so $[\D\RHom(F,F')] = [\D\RHom(h^*F, h^!F')] + [\D\RHom(j_{\cO_0}^*F, j_{\cO_0}^!F')]$.  The result then follows by induction on the number of orbits in the support of $F'$.
\end{proof}

\begin{lem}\label{lem:costalk-formula}
For any orbit $\cO \subset \cN$, we have
\[
[j_{\cO}^!\IC_\chi] = t^{-\dim \cO}\sum_{\{\psi \mid \cO_\psi = \cO\}} p_{\chi,\psi}(t^{-1})[L_\psi].
\]
\end{lem}
\begin{proof}
Using Corollary~\ref{cor:duality}, we have $j_{\cO}^!\IC_\chi \simeq \D j_{\cO}^* (\D \IC_\chi) \simeq \D (\IC_\chi|_{\cO})$.  We can therefore obtain a formula for for $[j_{\cO}^!\IC_\chi]$ by applying $\D$ to the right-hand side of~\eqref{eqn:p-defn}.  For any local system $L_\psi$ on $\cO$, we have $\D (L_\psi(-i)) \simeq L_\psi[2\dim \cO](\dim \cO + i)$, so the map $[\D(\cdot)]: K(\cO) \to K(\cO)$ sends $t^i[L_\psi] \mapsto t^{-\dim \cO - i}[L_\psi]$.  The result follows.
\end{proof}

\begin{lem}\label{lem:orbit-drhom}
We have
\[
[\D\RHom(j_{\cO}^*\IC_\chi, j_{\cO}^!\IC_\psi)] = \sum_{\{\phi,\phi' \mid \cO_\phi = \cO_{\phi'} = \cO\}} p_{\chi,\phi}(t)\lambda_{\phi,\phi'}(t)p_{\psi,\phi'}(t)[\bb_\pt].
\]
\end{lem}
\begin{proof}
Observe that $\D\RHom$ transforms Tate twists according to the formula $\D\RHom(F(n),F'(m)) \simeq \D\RHom(F,F')(n-m)$.  That means that the homomorphism $[\RHom(\cdot,\cdot)]: K(\cO) \times K(\cO) \to K(\pt)$ is $\Z[t^{1/2},t^{-1/2}]$-linear in the first variable, but \emph{antilinear} with respect to the involution $t^{1/2} \mapsto t^{-1/2}$ in the second variable.  Using~\eqref{eqn:p-defn} and Lemma~\ref{lem:costalk-formula}, we find that
\[
[\D\RHom(j_{\cO}^*\IC_\chi, j_{\cO}^!\IC_\psi)]
= 
t^{\dim \cO} \sum_{\!\!\!\!\!\{\phi,\phi' \mid \cO_\phi = \cO_{\phi'} = \cO\}\!\!\!\!\!} p_{\chi,\phi}(t) p_{\psi,\phi'}(t) [\D\RHom(L_\phi,L_{\phi'})].
\]
It suffices to show that $[\D\RHom(L_\phi,L_{\phi'})] = t^{-2\dim \cO}\lambda_{\phi,\phi'}(t)[\bb_{\pt}]$.  Using Corollary~\ref{cor:duality}, we have that $\RHom(L_\phi, L_{\phi'}) \simeq \RHom(\bb_{\cO}, L_\phi^* \otimes L_{\phi'}) \simeq R\Gamma(\cO, L_\phi \otimes L_{\phi'})$, where $L_\phi^*$ denotes the dual local system: $L_\phi^* = (\D L_\phi)[-2\dim \cO](-\dim \cO)$.  Therefore,
\[
\D\RHom(L_\phi,L_{\phi'}) \simeq R\Gamma_c(\cO, \D(L_\phi \otimes L_{\phi'}) \simeq R\Gamma(\cO, L_\phi \otimes L_{\phi'})[2\dim \cO](\dim \cO),
\]
so $[\D\RHom(L_\phi,L_{\phi'})] = t^{-\dim \cO}[R\Gamma_c(\cO, L_\phi \otimes L_{\phi'})]$, as desired.
\end{proof}

\begin{proof}[Proof of Theorem~\ref{thm:ls}]
The equation~\eqref{eqn:ls} follows from Lemmas~\ref{lem:orbit-sum} and~\ref{lem:orbit-drhom}, and the uniqueness assertion has been addressed in Remark~\ref{rmk:ls-alg}.
\end{proof}

\end{document}